\documentclass[12pt]{article}
\usepackage{amsfonts,amscd,amssymb}
\newtheorem{lemma}{Lemma}[section]
\newtheorem{theorem}[lemma]{Theorem}
\newtheorem{proposition}[lemma]{Proposition}
\newtheorem{corollary}[lemma]{Corollary}
\newtheorem{definition}[lemma]{Definition}
\newtheorem{remark}[lemma]{Remark}
\newtheorem{remarks}[lemma]{Remarks}

\newtheorem{examples}[lemma]{Examples}

\newenvironment{proof}{{\it Proof.}}{\hfill $ \square $ \vskip 4mm}

\newcommand{\thlabel}[1]{\label{th:#1}}
\newcommand{\thref}[1]{Theorem~\ref{th:#1}}
\newcommand{\selabel}[1]{\label{se:#1}}

\newcommand{\prlabel}[1]{\label{pr:#1}}
\newcommand{\prref}[1]{Proposition~\ref{pr:#1}}
\newcommand{\colabel}[1]{\label{co:#1}}
\newcommand{\coref}[1]{Corollary~\ref{co:#1}}

\newcommand{\eqlabel}[1]{\label{eq:#1}}
\newcommand{\eqref}[1]{(\ref{eq:#1})}

\newcommand{\End}{\rm{End}\,}

\def\lan{\langle}
\def\ran{\rangle}

\def\text#1{\mbox{{\rm #1}}}

\def\ot{\otimes}

\def\doublerightleft#1#2{{\lower.2ex\vbox{
\hbox{${\smash{\mathop{\longrightarrow}\limits^{#1}}}$}\vspace*{-4mm}
\hbox{${\smash{\mathop{\longleftarrow}\limits_{#2}}}$}}}}

\begin{document}
\title{Heisenberg double, pentagon equation,
structure and classification of finite dimensional Hopf algebras}
\author{G. Militaru\\
University of Bucharest\\ Faculty of Mathematics\\ Str. Academiei 14
\\RO-70109 Bucharest 1, Romania\\
email: gmilit@al.math.unibuc.ro }
\date{}
\maketitle
\begin{abstract}
\noindent
The study at the lavel of algebras of the pentagon equation
$R^{12}R^{13}R^{23}=R^{23}R^{12}$ leds
to the Structure and Classification theorem for finite
dimensional Hopf algebras.
We shall prove that $L$ is a finite dimensional Hopf algebra
if and only if there exists an invertible matrix
$R=(A_{ij})_{i,j=1,n} \in {\cal M}_{n^2}(k)\cong
{\cal M}_n (k)\ot {\cal M}_n (k)$,
solution of the pentagon equation
$\sum_{j=1}^n A_{ij}\ot A_{jp}=R(A_{ip}\ot I_n)R^{-1}$,
such that $L\cong P({\cal M}_n(k), R)=P(n, R)$.
Finally, we prove the Classification Theorem:
there exists a one to one correspondence between the
set of types of $n$-dimensional Hopf algebras and
the set of orbits of the action
$GL_n(k)\times ({\cal M}_n(k)\ot {\cal M}_n(k)) \to
{\cal M}_n(k)\ot {\cal M}_n(k)$, $(u, R)\to (u\ot u)R(u\ot u)^{-1}$,
the representatives of which are invertible solutions of
length $n$ for the pentagon equation.
\end{abstract}

\section*{Introduction}
The theory of Hopf algebras proved to be a unifying framework
for problems arising from different fields like: groups theory,
Galois theory, affine algebraic groups, Lie algebras,
local compact groups, operator algebras and quantum mechanics.
Recently, using the classification of finite groups as a model,
first steps were taken towards the classification of certain
types of finite dimensional Hopf algebras. The techniques
used are extremely diverse and the theory of Hopf algebras
classification is far from being completed (see the review
paper \cite{Mont} and its list of references).

In this paper we shall study, at the lavel of algebras, the
pentagon equation (also called the fusion equation)
$R^{12}R^{13}R^{23}=R^{23}R^{12}$, untill now analyzed mainly
in the theory of $C^*$-algebras (see \cite{BaajS93} and the
references indicated here).
Let $G$ be a locally compact group and {\sl dg} a right Haar
measure on $G$. Then the operator $V_G(\xi)(s,t)=\xi (st,t)$ is a
solution of the pentagon equation. Let $H$ be a separable Hilbert space
and $R\in {\cal L}(H\ot H)$ be a unitary solution of the
pentagon equation. If $R$ is commutative, then $R$ is equivalent
to $V_G\ot {\rm Id}_{K\ot K}$, for a locally compact group $G$
and a Hilbert space $K$ \cite[Theorem 2.2]{BaajS93}. If
$H$ is finite dimensional, then
\cite[Theorem 4.7 and Theorem 4.10]{BaajS93}
proves that $R$ is equivalent
to $V_S\ot {\rm Id}_{K\ot K}$ where, this time, $S$ is the
Woronowicz $C^*$-algebra associated to $R$. R. M. Kashaev shows
in \cite{Kashaev} that
the operator $V_G$ is in fact the "coordinate" representation of
the canonical element associated to the Heisenberg double of
the group algebra $kG$.
This remark will play an important role in proving
the structure of $\underline{Pent}$, the category of
pentagon objects: the objects of it are pairs $(A, R)$
where $A$ is a finite dimensional algebra and
$R\in A\ot A$ is an invertible solution of the pentagon equation
$R^{12}R^{13}R^{23}=R^{23}R^{12}$ in $A\ot A\ot A$.
Let $\underline{Hopf}$ be the category of finite dimensional Hopf
algebras, $L\in \underline{Hopf}$ and ${\cal H}(L)=L\# L^*$
the Heisenberg double of $L$. ${\cal H}(L)$ is an algebra
(we shall prove, using Morita theory, that it is isomorphic to
the matrix algebra ${\cal M}_{{\rm dim }(L)}(k)$) the
representations of which are the classic Hopf modules
${}^L{\cal M}_L$. \cite[Theorem 5.2]{DaeleK} and
\cite[Theorem 1]{Kashaev} proves that the canonical element
${\cal R}\in {\cal H}(L)\ot {\cal H}(L)$ of the Heisenberg double
is a solution of the pentagon equation.
Using a similar construction to the one of P. S. Baaj and
G. Skandalis (\cite{BaajS93}), we shall prove that there
exists a functor
$$
P: \underline{Pent}\to \underline{Hopf}, \quad
(A, R)\to P(A, R)
$$
such that $L\cong P({\cal H}(L), {\cal R})$ for any
$L\in \underline{Hopf}$. The structure of the solutions for
the pentagon equation shows that, up to an algebra map, these
solutions are canonical elements of varoius Heisebergs:
i.e. $(A, R)\in \underline{Pent}$
if and only if there exists $L\in \underline{Hopf}$ and
an algebra map $F: {\cal H}(L)\to A$ such that
$R=(F\ot F)({\cal R})$, where
${\cal R}\in {\cal H}(L)\ot {\cal H}(L)$ is the canonical
element of the Heisenberg double ${\cal H}(L)$.
The role of the Hopf modules category in solving the pentagon
equation was evidenced before in \cite[Proposition 1.3]{Street98}
and independently in \cite{M3}, where, in the equivalent context
of the Hopf equation, a FRT type theorem was proved.
Using the Fundamental Theorem for Hopf modules, a Lagrange type
theorem is proven: for $(A, R)\in \underline{Pent}$,
$A$ is free as a right $P(A, R)$-module and
${\rm dim}(A)={\rm dim}(P(A,R)){\rm dim}(A^{R,r})$,
where $A^{R,r}$ is the subspace of what we called
the right $R$-coinvariants of $A$.
The most important consequence of the results mentionated above
is the fact that they led us to the Structure and the Classification
theorem for finite dimensional Hopf algebras.
\thref{structurefinite} proves that
$L$ is a finite dimensional Hopf algebra if and only if there exists a
invertible matrix
$R \in {\cal M}_n (k)\ot {\cal M}_n (k)\cong {\cal M}_{n^2}(k)$,
solution of the pentagon equation, such that
$L\cong P({\cal M}_n (k), R)=P(n, k)$. The actual description,
in terms of $R$, of the Hopf algebra structure of $P(n, R)$
is obviuosly included: $P(n, R)$ is the subalgebra of
${\cal M}_n (k)$ of the left coefficients of $R$ with the
comultiplication given by
$\Delta (x)=R^{-1}(I_n\ot x)R$, for all $x\in P(n,R)$.
By writing $R$ as
$R=(A_{ij})_{i,j=1,n} \in {\cal M}_{n^2}(k)\cong
{\cal M}_n (k)\ot {\cal M}_n (k)$,
for some matrices $A_{ij}\in {\cal M}_n(k)$,
we get a much more elegant form of the pentagon equation for $R$
$$
\sum_{j=1}^n A_{ij}\ot A_{jp}=
R(A_{ip}\ot I_n)R^{-1}$$
for all $i$, $p=1,\cdots n$.\footnote{If $A$ and
$B\in {\cal M}_n(k)$, $A\ot B\in {\cal M}_{n^2}(k)$
is the Kronecker product of $A$ and $B$.}
\thref{thclas} is the Classification Theorem:
there exists a one to one correspondence between the
set of types of $n$-dimensional Hopf algebras and
the set of orbits of the action
$$
GL_n(k)\times ({\cal M}_n(k)\ot {\cal M}_n(k)) \to
{\cal M}_n(k)\ot {\cal M}_n(k), \quad
(u, R)\to (u\ot u)R(u\ot u)^{-1},
$$
the representatives of which are invertible solutions of
length $n$ for the pentagon equation. In this way, the theory
of finite dimensional Hopf algebras classification is reduced to
a new Jordan type theory, which we called
{\sl restricted Jordan theory}, to be studied in more depth in another
paper.

\section{Preliminaries}\selabel{1}
Throughout this paper, $k$ will be a commutative field. Unless
specified otherwise, all vector spaces, algebras, Hopf algebras,
tensor products and homomorphisms are over $k$. Let $V$ be a finite
dimensional vector space and
$\varphi: V\ot V^*\to \End(V)$,
$\varphi (v\ot v^*)(w)=\lan v^*, w\ran v$
be the canonical isomorphism. We recall that the element
$R:=\varphi^{-1}({\rm Id}_V)$ is called the
{\sl canonical element} of $V\ot V^*$. Of course
$$
R=\sum_{i=1}^n e_i\ot e_i^*
$$
where $\{e_i, e_i^* \mid i=1, \cdots, n\}$ is a dual basis and
$R$ is independent of the choice of the dual basis.

For a $k$-algebra $A$, ${{\cal M}}_A$ (resp. ${}_A{\cal M}$) will be
the category of right (resp.
left) $A$-modules and $A$-linear maps. ${\cal M}_n(k)$ will be the
$k$-algebra of $n\times n$-matrices and $GL_n(k)$ will be the
group of invertible $n\times n$-matrices. We shall denote by
$(e_{ij})$ the canonical basis of ${\cal M}_n(k)$:
$e_{ij}$ is the matrix having
$1$ in the $(i,j)$-position and $0$ elsewhere.
An $n$-dimensional $k$-algebra $A$ is viewed
as a subalgebra of ${\cal M}_n(k)$, via:
$A\subset \End (A)\cong {\cal M}_n(k)$.
If $A=(a_{ij})$, $B=(b_{ij}) \in {\cal M}_n(k)$ then, via
the canonical isomorphism
${\cal M}_n(k)\ot {\cal M}_n(k)\cong {\cal M}_{n^2}(k)$,
$A\ot B$ viewed as a matrix of ${\cal M}_{n^2}(k)$ is given
by the Kronecker product:
\begin{equation}\eqlabel{AoriB}
A\ot B=
\left(
\begin{array}{ccc}
a_{11}B&\cdots &a_{1n}B\\
 \cdot &\cdots &\cdot \\
a_{n1}B&\cdots &a_{nn}B
\end{array}
\right)
\end{equation}

Let $A$ be a $k$-algebra and $0\neq R=\sum R^1\ot R^2\in A\ot A$.
Let $R_{(l)}$, $R_{(r)}$ be the subspaces of left, respectively
right coefficients of $R$, i.e.
$$
R_{(l)}=\{\sum \lan a^*, R^2\ran R^1 \mid a^*\in A^* \}, \quad
R_{(r)}=\{\sum \lan a^*, R^1\ran R^2 \mid a^*\in A^* \}
$$
Let us assume that $R=\sum_{i=1}^m a_i\ot b_i$, where $m$ is
as small as possible. Then $m$ is called the {\sl length} of
$R$ and will be denoted $l(R)=m$. From the choice of $m$,
the sets $\{a_i\mid i=1,\cdots, m\}$, respectively
$\{b_i\mid i=1,\cdots, m\}$ are linear independent in $A$
and hence bases of $R_{(l)}$, respectively $R_{(r)}$. In particular,
dim$(R_{(l)})$=dim$(R_{(r)})$=$l(R)$.
Two elements $R$ and $S\in A\ot A$ are called {\sl equivalent}
(we shall write $R\sim S$) if there exists
$u\in U(A)$ an invertible element of $A$ such that
$S={}^uR:= (u\ot u)R(u\ot u)^{-1}$. If $R\sim S$ then
$l(R)$=$l(S)$. In particular,
if $\{ a_i \mid i=1, \cdots, m\}$ is a basis of $R_{(l)}$, then
$\{ ua_iu^{-1} \mid i=1, \cdots, m\}$ is a basis of
${}^uR_{(l)}$.

For coproducts and comodules we use Sweedler's notation
with suppressed summation sign: if $L$ is a
Hopf algebra, then for all $l\in L$ we write
$\Delta (l)=l_{(1)}\ot l_{(2)}$,
$(\Delta \ot {\rm Id})\Delta (l)=({\rm Id} \ot \Delta)\Delta (l)=
l_{(1)}\ot l_{(2)}\ot l_{(3)}$. If ($M,\rho_M$) is a left
$L$-comodule, then we write
$\rho_M(m)=m_{<-1>}\otimes m_{<0>} \in L\ot M$, for all
$m\in M$. ${}^L{{\cal M}}$ (resp. ${\cal M}^L$)
will be the category of left (resp. right) $L$-comodules and
$L$-colinear maps.

\subsection{The Heisenberg double and Hopf modules}
Let $L$ be a Hopf algebra with a bijective antipode.
We recall that a right-left $L$-Hopf module is a
threetuple $(M, \cdot, \rho)$, where
$(M, \cdot)$ is a right $L$-module, $(M, \rho)$ is a
left $L$-comodule such that the following compatibility
relation holds
\begin{equation}\eqlabel{comp}\eqlabel{1.1a}
\rho(m\cdot l)=\rho(m)\Delta(l)
\end{equation}
for all $l\in L$, $m\in M$. ${}^L{\cal M}_L$ will be the
category of right-left $L$-Hopf modules and $L$-linear, $L$-colinear
maps. If $V$ is a vector space, then
$V\ot L$ is a $L$-Hopf module via $(v\ot l)\cdot g:= v\ot lg$,
and $v\ot l\to l_{(1)}\ot v\ot l_{(2)}$, for all $v\in V$,
$l$, $g\in L$. Thus we have a functor
$-\ot L: {\cal M}_k \to {}^L{\cal M}_L$, $V\to V\ot L$.
Let
$(-)^{{\rm co}(L)}: {}^L{\cal M}_L\to {\cal M}_k$,
$M\to M^{{\rm co}(L)}=\{m\in M \mid \rho(m)=1_L \ot m\}$
be the functor that associates to each $M$ its coinvariants.
The Fundamental Theorem for Hopf
modules \cite{Sweedler69} (in the right-left version) shows
that for any $L$-Hopf module $M$ the multiplication map
$\mu :M^{{\rm co}(L)} \ot L \to M$,
$\mu (m\ot l)=m\cdot l$
is an isomorphism of $L$-Hopf modules with the inverse
$\beta :M \to M^{{\rm co}(L)} \ot L$,
$\beta (m)=m_{<0>}\cdot S^{-1}(m_{<-1>})\ot m_{<-2>}$.
Hence, $(-\ot L)\circ (-)^{{\rm co}(L)}\cong {\rm Id}_{{}^L{\cal M}_L}$.
As we also have that
$(-)^{{\rm co}(L)} \circ (-\ot L)\cong {\rm Id}_{{\cal M}_k}$,
we obtain that the functors
$-\ot L: {\cal M}_k \to {}^L{\cal M}_L$,
$(-)^{{\rm co}(L)}: {}^L{\cal M}_L\to {\cal M}_k$
define an equivalence of categories.

The same is valid for all left-right types of $L$-Hopf modules:
${}_L{\cal M}^L$, ${\cal M}_L^L$ or ${}^L_L{\cal M}$. For instance,
let ${}_L{\cal M}^L$ be the category of left-right $L$-Hopf modules
whose objects are left $L$-modules and right $L$-comodules $M$
such that $\rho (l\cdot m) =\Delta(l)\rho(m)$
for all $l\in L$, $m\in M$. For this category the
inverse of the multiplication map $L\ot M^{{\rm co}(L)} \to M$, in the
Fundamental Theorem, is given by
$\alpha (m)=m_{<2>} \ot S^{-1}(m_{<1>})\cdot m_{<0>}$,
for all $m\in M$. To conclude, if the antipode of $L$ is bijective,
then we have the equivalences of categories
$$
{}^L{\cal M}_L\cong {\cal M}_k\cong {}_L{\cal M}^L.
$$
We recall that $L^*$ is a left
$L$-module algebra in the usual way
$\lan h\cdot g^*, h'\ran =\lan g^*, h'h\ran$
for all $h$, $h'\in L$, $g^*\in L^*$.
The {\sl Heisenberg double} of $L$ is by definition the smash
product
$${\cal H}(L)=L\# L^*,$$
i.e. ${\cal H}(L)=L\ot L^*$ as a vector space, and the multiplication
is given by
$$
(h\# h^*)(g\# g^*)=h_{(2)}g \# h^*(h_{(1)}\cdot g^*)
=h_{(2)}g \# h^* \lan g^*, ? h_{(1)} \ran
$$
for all $h$, $g\in L$, $h^*$, $g^*\in L^*$.
${\cal H}(L)$ is an associative algebra with the unit
$1_L\# \varepsilon_L$. The Heisenberg double
${\cal H}(L)=L\# L^*$ is in fact a particular case
(for the Doi-Koppinen datum $(H, A, C)=(L,L,L)$)
of the general smash product $A\# C^*$ introduced in \cite{Doi92}
in the right-right version.
The right-left version of it and the above description of the
Heisenberg double is given in \cite{CIMZ}.
The canonical isomorphism of vector spaces
$$
\varphi: {\cal H}(L) \to \End(L)={\cal M}_{{\rm dim }(L)}(k), \quad
\varphi (l\# l^*)(g)=\lan l^*, g\ran l
$$
is not an algebra map. However, we can prove the following

\begin{proposition}\prlabel{anmorita}
Let $L$ be a finite dimensional Hopf algebra. Then there exists
an algebra isomorphism
$${\cal H}(L)\cong {\cal M}_{{\rm dim }(L)}(k).$$
\end{proposition}

\begin{proof}
As $L$ is finite dimensional, the functor
$$
T: {}^L{\cal M}_L \to {\cal M}_{{\cal H}(L)}, \quad
T(M)=M
$$
where the right ${\cal H}(L)$-action on $M$ is given by
$$
m\bullet (l\#l^*)=\lan l^*, m_{<-1>}\ran m_{<0>}\cdot l
$$
is an equivalence of categories (\cite{Doi92}).
As the antipode of $L$ is bijective (\cite{Sweedler69})
we have the following equivalences of categories
$$
{\cal M}_{{\cal H}(L)}\cong {}^L{\cal M}_L\cong {\cal M}_k
$$
i.e. ${\cal H}(L)$ is Morita equivalent to $k$. It follows from
the Morita theory
that there exists an algebra isomorphism
${\cal H}(L)\cong {\cal M}_n(k)$. Taking into account that
${\rm dim}({\cal H}(L))={\rm dim}(L)^2$, we obtain that $n={\rm dim}(L)$.
\end{proof}

\begin{remark}
\rm
The $k$-linear maps
$$
i_L :L\to {\cal H}(L), \;
i_L (l)=l\# \varepsilon_L \quad {\rm and} \quad
i_{L^*}: L^*\to {\cal H}(L), \;
i_{L^*} (l^*)=1_L\# l^*
$$
for all $l\in L$, $l^*\in L^*$ are injective algebra maps.
Beeing a smash product, the Heisenberg double ${\cal H}(L)$
satisfies the following universal property: given an algebra $A$
and the algebra maps $u: L\to A$, $v: L^*\to A$ such that
\begin{equation}\eqlabel{heisuniv}
u(l)v(l^*)=v(l_{(1)}\cdot l^*)u(l_{(2)})
\end{equation}
there exists a unique algebra
map $F: {\cal H}(L)\to A$ (given by
$F(l\#l^*):= v(l^*)u(l)$, for all $l\in L$,
$l^*\in L^*$) such that
\begin{equation}\eqlabel{dukdia}
F\circ i_L=u, \quad F\circ i_{L^*}=v
\end{equation}
\end{remark}

\subsection{The pentagon equation}\selabel{3}

\begin{definition}
Let $A$ be an algebra and $R=\sum R^1\ot R^2\in A\ot A$ be an
invertible element.
We shall say that $R$ is a solution of the
pentagon equation if
\begin{equation}\eqlabel{pentagon}
R^{12}R^{13}R^{23}=R^{23}R^{12}
\end{equation}
in $A\ot A\ot A$. $R$ is called unitary if
$1_A \in R_{(l)} \cap R_{(r)}$, i.e. there exist
$\lambda$, $\omega\in A^*$ such that
\begin{equation}\eqlabel{unitary}
\sum \lan \lambda, R^1\ran R^2=\sum \lan \omega, R^2\ran R^1=1_A
\end{equation}
\end{definition}

$\underline{Pent}$ will be the category of the (finite)
pentagon objects: the objects of it are pairs $(A, R)$,
where $A$ is a finite dimensional algebra and
$R\in A\ot A$ is an invertible solution of the pentagon equation.
A morphism $f: (A, R)\to (B, T)$ between two pentagon objects
$(A, R)$ and $(B, T)$ is an algebra map $f: A\to B$ such that
$(f\ot f)(R)=T$. $\underline{Pent}$ is a monoidal category
under the product $(A, R)\ot (B, T):=(A\ot B, R^{13}T^{24})$.

Let $R\in A\ot A$ be an invertible solution of the pentagon equation.
The subspaces of $A$
$$
A^{R,l}=\{a \in A \mid R (a\ot 1_A)=a\ot 1_A \} \;\;{\rm and}\;\;
A^{R,r}=\{a \in A \mid (1_A\ot a)R=1_A\ot a \}
$$
are called the spaces of left, respectively right,
{\sl $R$-coinvariants} of $A$.
On the other hand, the subspaces
$$
{}^{R, l}A=\{a \in A \mid R (1_A\ot a)=1_A\ot a \} \;\;{\rm and}\;\;
{}^{R, r}A=\{a \in A \mid (a\ot 1_A)R=a\ot 1_A \}
$$
are called the spaces of left, respectively right,
{\sl $R$-invariants} of $A$.

\begin{remarks}\label{2rem}
\rm
1. $R$ is a solution for the pentagon
equation if and only if $T=R^{-1}$ is a solution
of what we called in \cite{M3} the {\sl Hopf equation}:
$T^{12}T^{23}=T^{23}T^{13}T^{12}$.

2. Let $R\in A\ot A$ be an invertible solution of the
pentagon equation. If the subspaces $A^{R,l}$ and
$A^{R,r}$ are non-zero then $R$ is unitary. Indeed,
let $a\neq 0$ be an element of $A^{R,l}$ and
$R=\sum_{i=1}^m a_i\ot b_i\in A\ot A$, where $m=l(R)$. Then
$\sum_{i=1}^m a_i a \ot b_i =a\ot 1_A$ and hence,
$1_A=\sum_{i=1}^m \lan a^*, a_ia\ran b_i \in R_{(r)}$,
where $a^*\in A^*$ such that $\lan a^*, a\ran =1$. Similarly,
if $A^{R,r}\neq 0$ then $1_A \in R_{(l)}$.
It follows from \coref{unifd} that
an invertible solution of the pentagon equation is unitary
if and only if the spaces $A^{R,l}$ and $A^{R,r}$ are nonzero.
On the other hand, \coref{unifd} proves that if $A$
is finite dimensional, any invertible solution of the pentagon
equation is unitary.

3. In some applications, the algebra $A$ is of the form $A=\End_k(M)$,
where $M$ is a vector space over $k$. Then we can view $R$ as
an element of $\End_k(M\ot M)$, using the natural embedding
$\End_k(M)\ot \End_k(M)\subseteq \End_k(M\ot M)$,
and the pentagon equation as an equation in
$\End_k(M\ot M\ot M)$. It was proved in \cite[Proposition 7.1]{Davydov}
that there exists a bijective correspondence between all structures
of monoidal categories on the category of vector spaces over $k$
and all bijective solutions $R\in \End_k(M\ot M)$
of the pentagon equation on various $k$-vector spaces $M$.
\end{remarks}

\begin{examples}
\rm
1. Let $R\in A\ot A$ be an invertible solution of the pentagon equation
and $f: A\to B$ an algebra map. Then $(f\ot f)(R) \in B\ot B$
is an invertible solution of the pentagon equation.

2. If $R\in A\ot A$ is a unitary invertible solution of the
pentagon equation and $u\in U(A)$ is an invertible element of $A$,
then ${}^uR= (u\ot u)R(u\ot u)^{-1}$ is a unitary invertible
solution of the pentagon equation.

3. Let $A$ be an algebra $a$, $b\in A$ and $R\in A\ot A$ given by
$$
R=1_A\ot 1_A + a\ot b.
$$
Then $R$ is a solution of the pentagon equation if and only if
\begin{equation}\eqlabel{gi2nilsol}
a\ot (ab-ba-1)\ot b=
a\ot a\ot b^2 + a^2\ot b\ot b + a^2\ot ba\ot b^2.
\end{equation}
Now, if we take $a$ and $b$ such that
\begin{equation}\eqlabel{nilsol1}
a^2=b^2=0 \quad {\rm and} \quad ab-ba=1
\end{equation}
or
\begin{equation}\eqlabel{nilsol2}
a^2=0, \quad b^2=b \quad {\rm and} \quad
ab-ba=a+1
\end{equation}
then \eqref{gi2nilsol} is satisfied and $R$ is an invertible solution
of the pentagon equation.

We shall give some specific examples of these types in the case
that $A={\cal M}_n(k)$. First, in both cases,
it follows from the last equation that char$(k)|n$.

Assume now that char$(k)=2$ and let $n=2q$,
where $q$ is a positive integer. Then
$$
a=e_{12} + e_{34} +\cdots + e_{2q-1, 2q}, \quad
b=e_{21} + e_{43} +\cdots + e_{2q, 2q-1}
$$
is a solution of \eqref{nilsol1}, and hence
$R=I_n \ot I_n + \sum_{i,j=1}^q e_{2i-1, 2i}\ot e_{2j, 2j-1}$
is an invertible solution of the pentagon equation.
On the other hand, for an arbitrary invertible matrix
$X\in {\cal M}_q(k)$, the matrices $a$,
$b\in {\cal M}_n(k)$ given by
$$
b=e_{11} + e_{22} +\cdots e_{qq}, \quad
a=
\left(
\begin{array}{cc}
I_q & X^{-1}\\
X   & I_q
\end{array}
\right)
$$
are solutions of \eqref{nilsol2} and hence
$R=R_X=I_n\ot I_n + \sum_{i=1}^q a\ot e_{ii}$
is an invertible solution of the pentagon equation.

4. Let $(M,\cdot, \rho)\in {}_L{\cal M}^L$ be a $L$-Hopf module
over a Hopf algebra $L$. Then the map
\begin{equation}\eqlabel{gi2rares}
R_{(M,\cdot, \rho)}: M\ot M\to M\ot M, \quad
R_{(M,\cdot, \rho)}(m\ot n)=m_{<0>}\ot m_{<1>}\cdot n
\end{equation}
is a bijective solution of the pentagon equation in
$\End_k(M\ot M\ot M)$ (\cite{M3}, \cite{Street98}).
Conversely, if $M$ is finite dimensional and
$R\in \End_k(M\ot M)$ is a solution of the
pentagon equation, then there exists a bialgebra $B(R)$
such that $(M,\cdot ,\rho) \in {}_{B(R)}{\cal M}^{B(R)}$
and $R=R_{(M,\cdot ,\rho )}$; this is \cite[Theorem 3.1]{M3}
taking into account the equivalence between the Hopf and
the pentagon equation: $R$ is a solution of the
pentagon equation in $\End_k(M\ot M\ot M)$
if and only if $\tau R\tau$ is a solution
of the Hopf equation, where $\tau\in \End_k(M\ot M)$ is the
flip map.
\end{examples}

\begin{proposition}\prlabel{p1}
Let $A$ be an algebra and $R\in A\ot A$ an invertible solution
of the pentagon equation. Consider the comultiplications
$\Delta_r$, $\Delta_l: A\to A\ot A$, given by
\begin{equation}\eqlabel{ec12}
\Delta_r (a) = R^{-1} (1_A\ot a) R=\sum U^1R^1\ot U^2 a R^2
\end{equation}
\begin{equation}\eqlabel{ec11}
\Delta_l (a) = R(a\ot 1_A)R^{-1}=\sum R^1 a U^1\ot R^2 U^2
\end{equation}
where $U=\sum U^1\ot U^2 =R^{-1}$. Then
$A_r=(A, \cdot, \Delta_r)$ and $A_l =(A, \cdot, \Delta_l)$
are bialgebras without counit.
\end{proposition}

\begin{proof}
It is obvious that $\Delta_r$, $\Delta_l$ are algebra maps.
For $a\in A$ we have
$$
({\rm Id}\ot \Delta_r)\Delta_r (a)=
(R^{23})^{-1}(R^{13})^{-1} (1_A\ot 1_A \ot a)R^{13}R^{23}
$$
and
$$
(\Delta_r \ot {\rm Id})\Delta_r (a)=
(R^{12})^{-1}(R^{23})^{-1} (1_A\ot 1_A \ot a)R^{23}R^{12}
$$
so $\Delta_r$ is coassociative if and only if
$$
R^{23}R^{12}(R^{23})^{-1}(R^{13})^{-1} (1_A\ot 1_A \ot a)=
(1_A\ot 1_A \ot a) R^{23}R^{12}(R^{23})^{-1}(R^{13})^{-1}.
$$
Using the pentagon equation \eqref{pentagon}, we find that this
is equivalent to
$$
R^{12}(1_A\ot 1_A \ot a)=(1_A\ot 1_A \ot a) R^{12}
$$
and this equality holds for any $a\in A$. In a similar way we
can prove that $\Delta_l$ is also coassociative.
\end{proof}

It follows from \prref{p1} that we can put two different algebra
structures (without unit) on the dual $A^*$: the multiplications are
the convolutions $*_l$ and $*_r$ which are the dual maps of
$\Delta_l$ and  $\Delta_r$, i.e.
\begin{equation}\eqlabel{ec50}
\lan \omega *_l \omega', a\ran =\sum
\lan \omega, R^1 a U^1\ran \lan \omega', R^2 U^2\ran
\end{equation}
\begin{equation}
\lan \omega *_r \omega', a\ran =\sum
\lan \omega, U^1 R^1\ran \lan \omega', U^2 a R^2 \ran
\end{equation}
for all $\omega$, $\omega'\in A^*$, $a\in A$.

Aside from the invertibility  condition for
${\cal R}$,
the following theorem is \cite[Theorem 5.2]{DaeleK} and
\cite[Theorem 1]{Kashaev}. In \cite{DaeleK}, the Heisenberg double
does not appear explicitely, and in \cite{Kashaev} the
Heisenberg double is described in terms of structure constants,
and not as a smash product. For this reason, and for the
reader's convenience, we shall present the proof.

\begin{theorem}\thlabel{kdpen}
Let $L$ be a finite dimensional Hopf algebra and
$\{e_i, e_i^* \mid i=1,\cdots, n \}$ a dual basis of $L$.
Then the canonical element
$$
{\cal R}=\sum_i (e_i\# \varepsilon)\ot (1\# e_i^*)
\in {\cal H}(L)\ot {\cal H}(L)
$$
is an invertible solution of the pentagon equation in
${\cal H}(L)\ot {\cal H}(L)\ot {\cal H}(L)$.

Consequently, if $A$ is an algebra and
$f:{\cal H}(L) \to A$ an algebra map, then
$(f\ot f)({\cal R})$ is an invertible solution of
the pentagon equation in $A\ot A\ot A$.
\end{theorem}

\begin{proof}
Taking into account the multiplication rule of ${\cal H}(L)$ we have
$$
{\cal R}^{23}{\cal R}^{12}=
\sum_{i,j} (e_j\# \varepsilon)\ot
(e_{i_{(2)}}\# e_{i_{(1)}}\cdot e_j^*)\ot (1\# e_i^*)
$$
and
$$
{\cal R}^{12}{\cal R}^{13}{\cal R}^{23}=
\sum_{a,b,c} (e_a e_b\# \varepsilon)\ot (e_c\# e_a^*)\ot
(1\# e_b^*e_c^*)
$$
so we have to prove the equality
\begin{equation}\eqlabel{tre3pa}
\sum_{i,j} e_j\ot e_{i_{(2)}} \ot e_{i_{(1)}}\cdot e_j^* \ot e_i^*=
\sum_{a,b,c} e_ae_b \ot e_c\ot e_a^* \ot e_b^*e_c^*
\end{equation}
in $L\ot L\ot L^*\ot L^*$.
Let's fix the indices $x$, $y$, $z$ and $t=\{1,\cdots, n\}$ and
evaluate \eqref{tre3pa} at $e_x^*\ot e_y^*\ot e_z\ot e_t$.
\eqref{tre3pa} is then equivalent to
$$
\lan e_y^*, e_{t_{(2)}}\ran \lan e_x^*, e_z e_{t_{(1)}}\ran =
\sum_b \lan e_x^*, e_z e_b\ran \lan e_b^* e_y^*, e_t\ran,
$$
which can be obtained by applying the definition of the
convolution product and the dual basis formula
in the right hand side. We shall prove now that
$$
U=\sum_i (S(e_i)\# \varepsilon)\ot (1\# e_i^*)
$$
is the inverse of $R$, where $S$ is the antipode of $L$.
Since ${\cal H}(L)\ot {\cal H}(L)$ is isomorphic to
${\cal M}_{n^2}(k)$, it is enough to prove that
$RU=1\ot 1$. As
$RU=\sum_{i,j} (e_i S(e_j)\# \varepsilon)\ot (1\# e_i^*e_j^*)$,
we have to prove the formula
$$
\sum_{i,j} e_i S(e_j) \ot e_i^*e_j^*= 1\ot \varepsilon
$$
which holds, as for indices $x$, $y=1,\cdots, n$ we have
\begin{eqnarray*}
\sum_{i,j} \lan e_x^*, e_i S(e_j)\ran \lan e_i^*e_j^*, e_y\ran &=&
\sum_{i,j} \lan e_x^*, e_i S(e_j)\ran
\lan e_i^*, e_{y_{(1)}}\ran \lan e_j^*, e_{y_{(2)}}\ran \\
&=&\sum_{i,j} \lan e_x^*, e_i \lan e_i^*, e_{y_{(1)}}\ran
S(e_j \lan e_j^*, e_{y_{(2)}}\ran)\ran \\
&=& \lan e_x^*, e_{y_{(1)}}
S(e_{y_{(2)}}) \ran
= \lan e_x^*, 1\ran \lan \varepsilon, e_y\ran
\end{eqnarray*}
\end{proof}

\section{The solutions of the Pentagon equation}\selabel{4}
We shall set some notations that will remain
valid during the rest of the paper: $R=\sum R^1\ot R^2 \in A\ot A$
will be an invertible solution of the pentagon equation, where $A$ is
an algebra.
Let $R_{(l)}$, $R_{(r)}$ be the subspaces of left, respectively
right coefficients of $R$. We shall denote them as follows
$$
P=P(A,R):=R_{(l)}, \quad H=H(A,R):=R_{(r)}.
$$
Assume now that $R=\sum_{i=1}^m a_i\ot b_i$, where $m=l(R)$.
The fact that $m$ is minimal
implies that the sets $\{a_i\mid i=1,\cdots, m\}$ and
$\{b_i\mid i=1,\cdots, m\}$ are basis of $P=R_{(l)}$ and
$H=R_{(r)}$. In particular, dim$(P)$=dim$(H)=l(R)=m$.
Now consider $a_i^*\in P^*$ and $b_i^*\in H^*$ such that
$\{a_i, a_i^*\}$ and $\{b_i, b_i^*\}$ are dual bases of
$P$ and $H$:
$\lan a_i^*, a_j\ran =\delta_{ij}=\lan b_i^*, b_j\ran$.
Extend $a_i^* :P\to k$ and $b_i^* :H\to k$ to
respectively $\omega_i: A\to k$ and  $\lambda_i: A\to k$. We then
have
\begin{equation}\eqlabel{ec5}
\sum_{i=1}^m \lan \omega_k, a_i\ran b_i =b_k \quad {\rm and}
\quad
\sum_{j=1}^m a_j \lan \lambda_k, b_j\ran =a_k
\end{equation}
for all $k=1, \cdots, m$. We will use two different
notations for $R$:
$R=\sum_{i=1}^m a_i\ot b_i=\sum_{j=1}^m a_j\ot b_j$,
when we are interested in the basis elements of $P$ and $H$,
and the generic notation
$R=\sum R^1\ot R^2=\sum r^1\ot r^2=r$ and
$U=R^{-1}=\sum U^1\ot U^2 $.

The construction from part 1) of the next theorem is bassically
the one of P. S. Baaj and G. Skandalis \cite{BaajS93}, for
a unitary solution of the pentagon equation
$R\in {\cal L}(H\ot H)$, where $H$ is a separable Hilbert space.
Part 3) together with \thref{kdpen} give the structure of the
category $\underline{Pent}$.
Part 4) is a Lagrange type theorem and the
last part will play a role in the clasification of finite
dimensional Hopf algebras.

\begin{theorem}\thlabel{te1}
Let $A$ be an algebra, $R=\sum R^1\ot R^2 \in A\ot A$ an
unitary invertible solution of the pentagon equation and
$P=P(A,R)=R_{(l)}$, $H=H(A, R)=R_{(r)}$ the subspaces
of coefficients of $R$. Then:

1) $P$ and $H$ are subalgebras of $A$ and Hopf algebras
with the comultiplications given by the formulas
\begin{equation}
\Delta_P :P\to P\ot P, \quad
\Delta_P (x)=\Delta_r (x)=R^{-1}(1_A\ot x)R
\end{equation}
\begin{equation}\eqlabel{ec919}
\Delta_H :H\to H\ot H, \quad
\Delta_H (y)=\Delta_l (y)=R(y\ot 1_A)R^{-1}
\end{equation}
for all $x\in P$, $y\in H$. Furthermore,
the subalgebra $A_R$ of $A$ generated by $H$ and $P$
is finite dimensional and dim$(A_R)\leq l(R)^2$. \\
2) The $k$-linear map
$f: P^* \to H$, $f(p^*)=\sum \lan p^*, R^1\ran R^2$
is an isomorphism of Hopf algebras.\\
3) The $k$-linear map
$F: {\cal H}(P)\to A$,
$F(p\# p^*)=\sum \lan p^*, R^1\ran R^2 p$
is an algebra map and $R=(F\ot F)({\cal R})$,
where ${\cal R}\in {\cal H}(P)\ot {\cal H}(P)$ is
the canonical element associated to the Heisenberg double. \\
4) The multiplication on $A$ defines isomorphisms
$$
A^{R,r}\ot P \cong A \quad ({\rm resp.}\; H\ot A^{R,l}\cong A )
$$
of right $P$-modules (resp. left $H$-modules). In particular,
$A$ is free as a right $P$-module and as a left $H$-module and,
if $A$ is finite dimensional
$$
{\rm dim}(P)={\rm dim}(H)=\frac{{\rm dim}(A)}{{\rm dim}(A^{R,l})}=
\frac{{\rm dim}(A)}{{\rm dim}(A^{R,r})}.
$$
5) If $f: A\to B$ is an algebra isomorphism and
$S=(f\ot f)(R)$, then the Hopf algebras $P(A, R)$ and
$P(B, S)$ are isomorphic. Consequently, there exists an
Hopf algebra isomorphism
$P(A, R)\cong P(A, {}^uR)$, for any $u\in U(A)$.
\end{theorem}

\begin{proof}
1) We shall use the notations introduced above.
First we shall prove that $P$ (resp. $H$) are unitary
subalgebras in $A$ and subcoalgebras of $A_r=(A, \Delta_r)$
(resp. $A_l =(A, \Delta_l)$). This will follow from
the formulas:
\begin{equation}\eqlabel{ec7}
a_pa_q=
\sum_{j=1}^m \lan \lambda_p *_l \lambda_q, b_j\ran a_j \in P,
\quad
\Delta_r(a_p)=\sum_{i, j=1}^m \lan \lambda_p, b_ib_j\ran a_i\ot a_j
\in P\ot P
\end{equation}
and
\begin{equation}\eqlabel{ec9}
b_pb_q=\sum_{j=1}^m \lan \omega_p *_r \omega_q, a_j\ran b_j \in H,
\quad
\Delta_l(b_p)=\sum_{i, j=1}^m \lan \omega_p, a_ia_j\ran b_i\ot b_j
\in H\ot H
\end{equation}
for all $p$, $q=1, \cdots, m$. We prove \eqref{ec7},
(being similar, \eqref{ec9} is left to the reader).
\begin{eqnarray*}
\sum_{j=1}^m a_j\lan \lambda_p *_l \lambda_q, b_j\ran
&\stackrel{\eqref{ec50}}{=}&
\sum_{j=1}^m a_j \lan \lambda_p, R^1 b_j U^1\ran
\lan \lambda_q, R^2 U^2\ran \\
&=& ({\rm Id}\ot \lambda_p\ot \lambda_q)
(\sum_{j=1}^m a_j \ot R^1 b_j U^1 \ot R^2 U^2\ran )\\
&=& ({\rm Id}\ot \lambda_p\ot \lambda_q) (R^{23}R^{12}(R^{23})^{-1})
\stackrel{\eqref{pentagon}}{=}
({\rm Id}\ot \lambda_p\ot \lambda_q)(R^{12}R^{13})\\
&=& ({\rm Id}\ot \lambda_p\ot \lambda_q)
(\sum_{j,k=1}^m a_ja_k \ot b_j \ot b_k )\\
&=&\sum_{j,k=1}^m a_j \lan \lambda_p, b_j\ran a_k
\lan \lambda_q, b_k\ran
= a_pa_q
\end{eqnarray*}
i.e. $P$ is a subalgebra of $A$. On the other hand
\begin{eqnarray*}
\Delta_r (a_p)
&\stackrel{\eqref{ec5}}{=}&
\Delta_r (\sum_{j=1}^m a_j \lan \lambda_p, b_j\ran )
\stackrel{\eqref{ec12}}{=}
\sum_{j=1}^m U^1 R^1 \ot U^2 a_j R^2 \lan \lambda_p, b_j\ran \\
&=& ({\rm Id}\ot {\rm Id}\ot \lambda_p)
(\sum_{j=1}^m U^1 R^1 \ot U^2 a_j R^2 \ot b_j)\\
&=& ({\rm Id}\ot {\rm Id}\ot \lambda_p)
((R^{12})^{-1}R^{23}R^{12})
\stackrel{\eqref{pentagon}}{=}
({\rm Id}\ot {\rm Id}\ot \lambda_p)
(R^{12}R^{23})\\
&=& ({\rm Id}\ot {\rm Id}\ot \lambda_p)
(\sum_{i,j=1}^m a_i\ot a_j \ot b_i b_j)
= \sum_{i,j=1}^m a_i\ot a_j \lan \lambda_p, b_i b_j\ran
\end{eqnarray*}
i.e. $P$ is a subcoalgebra in $(A, \Delta_r)$.
A similar computation yields \eqref{ec9},
proving that $H$ is a subalgebra of $A$ and a subcoalgebra in
$(A, \Delta_l)$. As $R$ is unitary, $1_A\in P\cap H$, i.e.
$P$ and $H$ are unitary subalgebras of $A$.
Moreover, from the construction we have that
$R\in P\ot H\subset A\ot A$ and hence we can view
$U=R^{-1}\in P\ot H$.

We shall define now the counit and the antipode of the
Hopf algebras $P$ and $H$. They are given by the formulas:
\begin{equation} \eqlabel{ec90}
\varepsilon_P :P\to k, \quad
\varepsilon_P (a_k)=\lan b_k^*, 1_A \ran, \quad
S_P :P \to P, \quad
S_P (a_k)= \sum U^1 \lan b_k^*, U^2\ran
\end{equation}
and
\begin{equation}\eqlabel{ec91}
\varepsilon_H :H\to k, \quad
\varepsilon_H (b_k)=\lan a_k^*, 1_A \ran, \quad
S_H :H \to H, \quad
S_H (b_k)= \sum \lan a_k^* ,U^1 \ran U^2
\end{equation}
for all $k=1, \cdots, m$. We shall prove that $P$ is a
Hopf algebra; the fact that $H$ is a Hopf algebra
is proved in a similar way. First, we remark that, as $H$ is a
subalgebra of $A$, the comultiplication from
\eqref{ec7} can be rewritten as
\begin{equation}\eqlabel{ec81}
\Delta_r(a_p)=\sum_{i, j=1}^m \lan b_p^*, b_ib_j\ran a_i\ot a_j
\end{equation}
Now, for $p=1, \cdots, m$ we have
\begin{eqnarray*}
({\rm Id}\ot \varepsilon_P)\Delta_r (a_p)
&\stackrel{\eqref{ec81}}{=}&
\sum_{i,j=1}^m a_i \lan b_p^*, b_ib_j\ran \lan b_j^*, 1_A\ran
=
\sum_{i,j=1}^m a_i \lan b_p^*, b_ib_j \lan b_j^*, 1_A\ran \ran\\
&=&\sum_{i=1}^m a_i \lan b_p^*, b_i\ran =a_p
\end{eqnarray*}
i.e. $({\rm Id}\ot \varepsilon_P)\Delta_r ={\rm Id}$. A
similar computation shows that
$(\varepsilon_P \ot {\rm Id})\Delta_r ={\rm Id}$, and
$\varepsilon_P$ is a counit of $P$. $S_P$ is a right
convolution inverse of ${\rm Id}_P$ since
\begin{eqnarray*}
({\rm Id}\ot S_P)\Delta_r (a_p)
&\stackrel{\eqref{ec81}}{=}&
\sum_{i,j=1}^n a_i\lan b_p^*, b_ib_j\ran U^1 \lan b_j^*, U^2\ran
= \sum_{i,j=1}^n a_i U^1 \lan b_p^*, b_ib_j \lan b_j^*, U^2\ran\ran\\
&=&\sum_{i=1}^n a_i U^1 \lan b_p^*, b_i U^2\ran
=({\rm Id}\ot b_p^*) (RR^{-1})
= 1_A \lan b_p^*, 1_A\ran
= \varepsilon_P (a_p)1_A
\end{eqnarray*}
From the fact that $P$ is finite dimensional, it follows that
$S_P$ is an antipode of $P$.

We shall prove now that the subalgebra $A_R$ of
$A$ generated by $H$ and $P$ is finite dimensional.
We shall use the pentagon equation. We have:
$$
R^{12}R^{13}R^{23}=
\sum_{i,j,k=1}^m a_i a_j \ot b_ia_k \ot b_jb_k \quad
{\rm and} \quad
R^{23}R^{12}=\sum_{i,j=1}^m a_j \ot a_i b_j \ot b_i
$$
Applying $a_u^* \ot {\rm Id} \ot b_v^*$ to the
pentagon equation we obtain
$$
a_u b_v =\sum_{i,j,k=1}^m
\lan a_u^*, a_i a_j \ran \lan b_v^*, b_jb_k\ran b_ia_k
$$
for all $u$, $v=1, \cdots, m$. This formula gives that
dim$(A_R)\leq m^2$.

2) Let us prove now that $f: P^* \to H$ is an isomorphism of
Hopf algebras. $f$ is an isomorphism of vector spaces since
$f(a_j^*)=\sum_{i=1}^m \lan a_j^*, a_i\ran b_i=b_j$.
The formula \eqref{ec9} can be rewritten as
$$
b_pb_q=\sum_{j=1}^m \lan a_p^* *_r a_q^*, a_j\ran b_j
$$
which means that $f(a_p^*)f(a_q^*)=f(a_p^* *_r a_q^*)$
for all $p$, $q=1,\cdots, m$, i.e. $f$ is an algebra isomorphism.
Let us prove now that $f$ is also a coalgebra map. We recall the
definition of the comultiplication $\Delta_{P^*}$:
$$
\Delta_{P^*} (a_p^*)=\sum X^1\ot X^2 \in P^*\ot P^*, \;\;
{\rm if~and~only~if}\;\;
\lan a_p^*, xy\ran =\sum \lan X^1, x\ran \lan X^2, y\ran
$$
for all $x$, $y\in P$. It follows that
\begin{eqnarray*}
(f\ot f)\Delta_{P^*}(a_p^*)&=&
\sum f(X^1)\ot f(X^2)
=\sum \lan X^1, R^1\ran R^2 \ot \lan X^2, r^1\ran r^2\\
&=&\sum \lan a_p^*, R^1r^1\ran R^2\ot r^2
=\sum_{i,j=1}^m \lan a_p^*, a_ia_j\ran b_i\ot b_j\\
&\stackrel{\eqref{ec9}}{=}&
\Delta_H (b_p)
= (\Delta_H \circ f)(a_p^*)
\end{eqnarray*}
i.e. $f$ is also a coalgebra map. Hence, we have proved that
$f$ is an isomorphism of bialgebras and, as $P$ and $H$ are Hopf
algebras, it is also a isomorphism of Hopf algebras
(\cite{Sweedler69}).

3) We remark that
$$
F(a_i \# a_j^*)=\sum_{t=1}^m \lan a_j^*, a_t\ran b_ta_i=
b_ja_i
$$
for all $i$, $j=1,\cdots, m$. The fact that $F$ is an algebra
map can be proved directly by using this formula; another way to
proceed is to use the universal property of the Heisenberg double
${\cal H}(P)$ for the diagram \eqref{dukdia}, with $L=P$,
$u :P \to A$ is the usual inclusion and $v: P^*\to A$ is
the composition $v=f\circ j$, where $f:P^*\to H$ is the
isomorphism from part 2) and $j: H\to A$ is the usual inclusion.
We only have to prove that the compatibility
condition \eqref{heisuniv} holds, i.e.
$$
hv(g^*)=v(h_{(1)}\cdot g^*)h_{(2)}
$$
for any $h\in P$ and $g^*\in P^*$, which turns out to be
$$
\sum h \lan g^*, R^1\ran R^2=
\sum \lan g^*, R^1 h_{(1)}\ran R^2 h_{(2)}
$$
or, equivalently
$$
\sum R^1\ot hR^2=\sum R^1 h_{(1)}\ot R^2 h_{(2)}.
$$
This equation holds, as $\Delta_P (h)=R^{-1}(1_H\ot h)R$,
for any $h\in P$.

Now let ${\cal R}=\sum_{i=1}^m (a_i\# \varepsilon_P) \ot (1_A\# a_i^*)$
be the canonical element of ${\cal H}(P)\ot {\cal H}(P)$. Then
$$
(F\ot F)({\cal R})=\sum_{i,t=1}^m \lan \varepsilon_P, a_t\ran
b_t a_i\ot b_i=
\sum_{i,t=1}^m \lan b_t^*, 1_A\ran b_t a_i\ot b_i=
\sum_{i=1}^m a_i\ot b_i =R.
$$
4) Consider the map
$$
\rho=\rho_P: A\to P\ot A, \quad
\rho(a)=(1_A\ot a)R=\sum R^1\ot aR^2=\sum_{i=1}^m a_i \ot ab_i
$$
for all $a\in A$. We will show that
$(A, \cdot, \rho_P)\in {}^P{\cal M}_P$ is a right-left $P$-Hopf
module, where the structure of right $P$-module is simply the
multiplication $\cdot$ of $A$. Indeed, for $a\in A$ we have
\begin{eqnarray*}
({\rm Id}\ot \rho)\rho (a)&=&\sum R^1 \ot \rho(aR^2)
=\sum R^1\ot r^1 \ot a R^2r^2\\
&=&(1_A\ot 1_A \ot a)R^{13}R^{23}
\stackrel{\eqref{pentagon}}{=}
(1_A\ot 1_A \ot a) (R^{12})^{-1} R^{23}R^{12}\\
&=& \sum U^1r^1\ot U^2R^1r^2\ot aR^2
= \sum R^{-1}(1_A\ot R^1)R\ot aR^2\\
&=& \sum \Delta_P (R^1)\ot aR^2
= (\Delta_P \ot {\rm Id})\rho (a)
\end{eqnarray*}
and
$$
\sum_{i=1}^m \lan \varepsilon_P, a_i\ran ab_i=
\sum_{i=1}^m \lan b_i^*, 1_A\ran ab_i =a
$$
so $(A, \rho)$ is a left $P$-comodule. The
compatibility relation
$$
\rho(a)\Delta_P(a_i)=(1_A\ot a)RR^{-1}(1_A\ot a_i)R=
(1_A\ot aa_i)R=\rho(aa_i)
$$
holds for all $i=1, \cdots, m$ and $a\in A$. Hence,
$(A, \cdot, \rho_P)\in {}^P{\cal M}_P$ and the coinvariants
$$
A^{{\rm co}(P)}=\{a\in A\mid \rho(a)=1\ot a \}=A^{R,r}
$$
are the right $R$-coinvariants of $A$. From the
right-left version of the Fundamental Theorem of Hopf
modules it follows that the multiplication of $A$,
$$
\mu: A^{R,r}\ot P\to A, \quad \mu(a\ot x)=ax
$$
defines an isomorphism of $P$-Hopf modules and, in particular, of right
$P$-modules. We recall that $A^{R,r}\ot P$ is a right $P$-module via
$(a\ot x)\cdot y=a\ot xy$, for all $a\in A^{R,r}$, $x$, $y\in P$.
It follows that $A$ is free as a right $P$-module and, if
$A$ is finite dimensional,
$$
{\rm dim}(A)={\rm dim }(P){\rm dim}(A^{R,r}).
$$
In a similar way we can show that
$(A, \cdot, \rho_H) \in {}_H{\cal M}^H$, where $\cdot$ is the
multiplication of $A$ and
$$
\rho_H :A \to A\ot H, \quad
\rho_H (a)=R(a\ot 1_A)=\sum R^1a\ot R^2
$$
for all $a\in A$. Moreover, $A^{{\rm co}(H)}=A^{R,l}$. If we apply
once again the Fundamental Theorem of Hopf modules (this time the
left-right version) we obtain the other part of the statement.

5) $S=(f\ot f)(R)=\sum_{i=1}^m f(a_i) \ot f(b_i)$, $l(S)=l(R)$
and $S^{-1}=\sum f(U^1)\ot f(U^2)$.
It follows that $\{f(a_i) \mid i=1, \cdots, m\}$ is a basis of
$P(B, S)$ and hence
the restriction of $f$ to $P(A, R)$ gives an algebra isomorphism
between $P(A, R)$ and $P(B, S)$ that is also a coalgebra map since
\begin{eqnarray*}
(f\ot f)\Delta_{P(A, R)}(a_i)&=&
\sum f(U^1) f(R^1) \ot f(U^2) f(a_i) f(R^2)\\
&=&S^{-1} (1\ot f(a_i))S
=\Delta_{P(B,S)} (f(a_i))
\end{eqnarray*}
for all $i=1,\cdots, m$. The last statement is obtain for
$B=A$ and $f: A\to A$, $f(x)=uxu^{-1}$ for all $x\in A$.
\end{proof}

\begin{corollary}\colabel{unifd}
Let $A$ be an algebra and $R\in A\ot A$ an invertible solution
of the pentagon equation. Then:

1) $R$ is unitary if and only if the subspaces
$A^{R,l}$ and $A^{R,r}$ are nonzero.\\
2) If $A$ is finite dimensional, then $R$ is unitary.
\end{corollary}

\begin{proof}
1. We have proved in Remark \ref{2rem} that if the spaces
$A^{R,l}$ and $A^{R,r}$ are nonzero, then $R$ is unitary.
Conversely, if $R$ is unitary, the isomorphisms
of vector spaces from part 4) of \thref{te1} show us that
$A^{R,l}$ and $A^{R,r}$ are nonzero.

2. Assume now that $A$ is finite dimensional.
From the construction of \thref{te1} we have that
$R\in P\ot H$. In particular, for any positive integer $t$
there exist scalars $\alpha_{ij} \in k$ such that
\begin{equation}\eqlabel{ec500}
R^t = \sum_{i,j=1}^m \alpha_{ij} a_i\ot b_j
\end{equation}
We shall prove now that $R$ is unitary, i.e. $1_A\in P\cap H$.
As $A$ is finite dimensional, $A$ can be embedded into a matrix algebra
$A\subset {\cal M}_n (k)$, where $n=$dim$(A)$. We consider
$$
R\in A\ot A\subset {\cal M}_n(k)\ot {\cal M}_n(k)\cong {\cal M}_{n^2}(k)
$$
As $R$ is invertible, using the Hamilton-Cayley theorem
in ${\cal M}_{n^2}(k)$, the identity matrix $I_{n^2}$ can be
represented as a linear combination of powers of $R$.
Hence, using \eqref{ec500}, we obtain in $A\ot A$ a linear combination
$1_A\ot 1_A =\sum_{i,j=1}^m \gamma_{ij} a_i\ot b_j$,
for some scalars $\gamma_{ij}\in k$. Therefore
$$
1_A=\sum_{i,j=1}^m a_i \lan 1_A^*, \gamma_{ij} b_j\ran =
\sum_{i,j=1}^m \lan 1_A^*, \gamma_{ij} a_i\ran b_j \in P\cap H
$$
i.e. $R$ is unitary.
\end{proof}

Let $R$ be a unitary invertible solution of the pentagon
equation. We have seen in \thref{te1} that the subalgebra
$A_R$ generated by the coefficient spaces $R_{(l)}$, $R_{(r)}$
is finite dimensional and, of course, we can
view $R\in A_R\ot A_R$. From this reason, in studying of the
unitary invertible solutions of the
pentagon equation, it is enough to focus on the case when the
that the algebra $A$ is finite dimensional. In this case, the
unitary condition of $R$ follows automatically.

Using \thref{kdpen} and \thref{te1} we obtain the following
Corollary, which is the algebraic version of
\cite[Theorem 4.7]{BaajS93}: the role of the operator $V_S$
is played by the canonical element of a Heisenberg double.

\begin{corollary}\colabel{ancorst}
$(A, R)\in \underline{Pent}$ if and only if there exists a
finite dimensional Hopf algebra $L$ and an algebra map
$F: {\cal H}(L)\to A$ such that $R=(F\ot F)({\cal R})$.
\end{corollary}

\begin{remark}
\rm
Let $(A, R)\in \underline{Pent}$. It follows from
\coref{ancorst} and \prref{anmorita} that there exists an algebra
map $F:{\cal M}_m(k)\to A$, where $m=l(R)$; $F$ is injective
since ${\cal M}_m(k)$ is a simple algebra.
Let $0\neq a_{ij}=F(e_{ij})\in A$, $i$, $j=1,\cdots, m$; then
$a_{ij} a_{kl}=\delta_{jk} a_{il}$ and $1_A =\sum_{i=1}^m a_{ii}$.
It follows from the Reconstruction Theorem of the matrix algebra
(\cite[Theorem 17.5]{Lam}) that there exists an algebra isomorphism
$$
A\cong {\cal M}_m(B), \;  {\rm where} \;
B=\{x\in A \; \mid \; xa_{ij}=a_{ij}x, \; \forall i, j=1, \cdots, m\}.
$$
Hence $A$ is a noncommutative algebra if $R$ is non-trivial
($l(R)>1$ or, equivalently, $R\neq 1_A\ot 1_A$). Furthermore,
${\rm dim}(A)=m^2 {\rm dim}(B)$ and hence, $l(R)^2|{\rm dim}(A)$.
\end{remark}

We are now going to describe the space of integrals of the Hopf
algebras $P=P(A, R)$ and $H=H(A, R)$, where
$(A,R) \in \underline{Pent}$. Since $A$ is free as a right
$P$-module and as a left $H$-module, there exists
$$
\pi_P :A \to P, \quad ({\rm respectiv}\quad  \pi_H :A\to H )
$$
a non-zero right $P$-linear (respectiv left $H$-linear) map.

\begin{proposition}
Let $(A,R) \in \underline{Pent}$. Then:

1) $\pi_P (a)$ is a right integral in $P(A,R)$ for any right
$R$-invariant $a\in {}^{R, r}A$.\\
2) $\pi_H (a)$ is a left integral in $H(A, R)$ for any left
$R$-invariant $a\in {}^{R, l}A$.
\end{proposition}

\begin{proof}
1. $a$ is a right $R$-invariant; hence,
$\sum_{i=1}^m a a_i\ot b_i =a\ot 1_A$.
If we apply $b_p^*$ to this equality we get that
$a a_p =\lan b_p^*, 1_A\ran a= \varepsilon_P (a_p)a$,
for any $p=1, \cdots, m$. As $(a_i)$ is a basis of $P$, we obtain
that $a x = \varepsilon_P (x)a$, for all $x\in P$. If we apply
the right $P$-module map $\pi_P$ to this, we obtain that
$\pi_P (a)$ is a right integral of $P$.\\
2. Left to the reader.
\end{proof}

\section{The structure and the classification of finite dimensional
Hopf algebras}\selabel{5}

We shall now prove the following Spliting Theorem:

\begin{theorem}\thlabel{te2}
Let $L$ be a finite dimensional Hopf algebra. Then
there exists an isomorphism of Hopf algebras
$$
L\cong P({\cal H}(L), {\cal R})
$$
where ${\cal R}$ is the canonical element of the
Heisenberg double ${\cal H}(L)$.
\end{theorem}

\begin{proof}
Let $\{e_i, e_i^* \mid i=1,\cdots, m \}$ be a basis of $L$ and
$$
{\cal R}=\sum_{i=1}^m (e_i\# \varepsilon_L) \ot (1_L\# e_i^*)
\in {\cal H}(L)\ot {\cal H}(L)
$$
the canonical element. We have to prove that the
Hopf algebra $P({\cal H}(L), {\cal R})$ extracted from part 1)
of \thref{te1} is isomorphic to $L$, with the initial
Hopf algebra structure. $i_L: L\to {\cal H}(L)$,
$i_L (l)=l\# \varepsilon_L $
is an injective algebra map. We identify
$$L\cong {\rm Im}(i_L)=L\# \varepsilon_L.$$
From the construction,
$P({\cal H}(L), {\cal R})$ is the subalgebra of
${\cal H}(L)$ having $(e_i\# \varepsilon_L)_{i=1, \cdots, m}$
as a basis; i.e. there exists an algebra isomorphism
$L\cong {\rm Im}(i_L)=P({\cal H}(L), {\cal R})$.
It remains to be proven that the coalgebra structure
(resp. the antipode) of $P({\cal H}(L), {\cal R})$
extracted from \thref{te1} is exactly the original
coalgebra structure (resp. the antipode) of $L$.
As the counit and the antipode of a Hopf algebra are uniquely
determined by the multiplication and the comultiplication, the only
thing left to be shown is the fact that, via the above
identification, $\Delta_P=\Delta_L$. This means that
$$
\Delta_L (e_i\#\varepsilon_L)=
{\cal R}^{-1} (1_{{\cal H}(L)}\ot e_i\#\varepsilon_L) {\cal R}
$$
or equivalently,
$$
{\cal R}\Delta_L (e_i\#\varepsilon_L)=
\Bigl( (1_L\# \varepsilon_L) \ot (e_i\# \varepsilon_L)\Bigl) {\cal R}.
$$
Now we compute
\begin{eqnarray*}
\Bigl( (1_L\# \varepsilon_L) \ot (e_i\# \varepsilon_L)\Bigl){\cal R}&=&
\Bigl( (1_L\# \varepsilon_L) \ot (e_i\# \varepsilon_L)\Bigl)
\Bigl( \sum_{j=1}^m (e_j\# \varepsilon_L) \ot (1_L\# e_j^*)\Bigl)\\
&=& \sum_{j=1}^m (e_j\# \varepsilon_L) \ot
(e_{i_{(2)}}\# e_{i_{(1)}} \cdot e_j^*)
\end{eqnarray*}
On the other hand
\begin{eqnarray*}
{\cal R}\Delta_L (e_i\# \varepsilon_L)&=&
\sum_{j=1}^m \Bigl( (e_j\# \varepsilon_L) \ot (1_L\# e_j^*)\Bigl)
\Bigl( (e_{i_{(1)}}\# \varepsilon_L)\ot
( e_{i_{(2)}}\# \varepsilon_L ) \Bigl)\\
&=&\sum_{j=1}^m (e_{j} e_{i_{(1)}}\# \varepsilon_L )\ot
(e_{i_{(2)}}\# e_j^*)
\end{eqnarray*}
Hence, we have to show the formula
\begin{equation}\eqlabel{an646}
\sum_{j=1}^m e_{j} e_{i_{(1)}}\ot e_{i_{(2)}} \ot e_j^*=
\sum_{j=1}^m e_j \ot e_{i_{(2)}} \ot e_{i_{(1)}} \cdot e_j^*
\end{equation}
Let us set the indices $a$, $b$, $k=1, \cdots, m$, and evaluate
\eqref{an646} at $e_a^* \ot e_b^*\ot e_k$. \eqref{an646}
is then equivalent to
$$
\lan e_a^*, e_ke_{i_{(1)}}\ran \lan e_b^*, e_{i_{(2)}}\ran =
\sum_{j=1}^m \lan e_a^*, e_j\ran \lan e_b^*, e_{i_{(2)}}\ran
\lan e_{i_{(1)}} \cdot e_j^*, e_k \ran
$$
and this is easily verified using the dual basis formula.
It follows that $\Delta_L=\Delta_P$ and
$L\cong P({\cal H}(L), {\cal R})$ as Hopf algebras.
\end{proof}

We have now arrived at the Structure Theorem for finite dimensional
Hopf algebras. Let $L$ be a finite dimensional Hopf algebra.
\prref{anmorita} proves that
there exists an algebra isomorphism
${\cal H}(L) \cong {\cal M}_n(k)$, where $n={\rm dim}(L)$.
Via this isomorphism the canonical element
${\cal R}\in {\cal H}(L)\ot {\cal H}(L)$ is viewed as
an element of ${\cal M}_n(k)\ot {\cal M}_n(k)$, or as a
matrix of ${\cal M}_{n^2}(k)$.
We shall now give the data which shows us how any finite dimensional
Hopf algebra is constructed.
Let $R\in {\cal M}_n(k)\ot {\cal M}_n(k)$ be
an invertible solution of the pentagon equation.
We write $R=\sum_{i=1}^m A_i \ot B_i$ where $m=l(R)$; then
the sets of matrices
$\{A_i \mid i=1, \cdots m\}$ and $\{B_i \mid i=1, \cdots m\}$
are linearly independent over $k$.
Let $\{B_i^* \mid i=1, \cdots m\}$ be the dual basis of
$\{B_i \mid i=1, \cdots m\}$ and
$U=R^{-1}=\sum U^1\ot U^1$.

The Hopf algebra $P({\cal M}_n(k), R)$ will be denoted
for simplicity by $P(n, R)$ and is described as follows:

$\bullet$ as an algebra, $P(n, R)$ is the subalgebra of
the $n\times n$-matrix algebra ${\cal M}_n(k)$ with
$\{A_i \mid i=1, \cdots m\}$ as a $k$-basis;

$\bullet$ the coalgebra structure and the antipode of
$P(n, R)$ are given by the following formulas:
\begin{equation}\eqlabel{craiov}
\Delta : P(n, R)\to P(n, R)\ot P(n, R),
\quad \Delta(A_p)=R^{-1}(I_n\ot A_p)R
\end{equation}
\begin{equation}
\varepsilon: P(n, R)\to k, \quad
\varepsilon (A_p)=\lan B_p^*, I_n\ran
\end{equation}
\begin{equation}
S:P(n, R) \to P(n, R), \quad
S(A_p)=\sum \lan B_p^*, U^2\ran U^1
\end{equation}
for all $p=1, \cdots, m$. The dual of $P(n, R)$ is the subalgebra
$H(n, R)=H({\cal M}_n(k), R)$ of the matrix algebra ${\cal M}_n(k)$
having $\{B_i \mid i=1, \cdots m\}$ as a $k$-basis and
the comultiplication given by
$$
\Delta_l :H(n, R)\to H(n, R)\ot H(n, R), \quad
\Delta_l (B_p)=R(B_p\ot I_n)R^{-1}
$$
for all $p=1, \cdots, m$.

Theorems \ref{th:te1} and \ref{th:te2} imply the
following Structure Theorem for finite
dimensional Hopf algebras.

\begin{theorem}\thlabel{structurefinite}
$L$ is a finite dimensional Hopf algebra if and only if
there exists a positive integer $n$ and an
invertible solution of the pentagon equation
$R\in {\cal M}_n(k)\ot {\cal M}_n(k)\cong {\cal M}_{n^2}(k)$
such that $L\cong P(n, R)$. Furthermore,
\begin{equation}\eqlabel{an32}
{\rm dim}(L)=\frac{n^2}{{\rm dim}({\cal M}_n(k)^{R,r})}=l(R).
\end{equation}
\end{theorem}

Let $n$ be a positive integer. The algebra isomorphism
${\cal H}(L)\cong {\cal M}_{{\rm dim}(L)}(k)$ and
\thref{te2} show us that an $n$-dimensional Hopf algebra
$L$ is isomorphic to a $P(n, R)$ for
$R\in {\cal M}_n(k)\ot {\cal M}_n(k)$
an invertible solution of the pentagon equation such that $l(R)=n$.
Overall, each solution of the pentagon equation belonging to
${\cal M}_n(k)\ot {\cal M}_n(k)$ determines a Hopf algebra, the
dimension of which is not necessarily $n$, but a divisor of
$n^2$ based on the Lagrange formula \eqref{an32}.

We are now going to prove the Classification Theorem for finite
dimensional Hopf algebras. Let $\underline{Pent}_n$ be the set
$$
\underline{Pent}_n=\{R\in {\cal M}_n(k)\ot {\cal M}_n(k) \;\mid\;
({\cal M}_n(k), R)\in \underline{Pent} \;{\rm and}\; l(R)=n \}.
$$
\begin{theorem}\thlabel{thclas}
Let $n$ be a positive integer. Then there exists a one to one
correspondence between the set of types of $n$-dimensional Hopf
algebras and the set of the orbits of the action
\begin{equation}\eqlabel{anactcon}
GL_n(k)\times \underline{Pent}_n \to \underline{Pent}_n, \quad
(u, R)\to (u\ot u)R(u\ot u)^{-1}.
\end{equation}
\end{theorem}

\begin{proof}
In part 5) of the \thref{te1} we proved that there exists a Hopf
algebra isomorphism $P(n, R)\cong P(n, {}^uR)$ for any
$R\in \underline{Pent}_n$ and
$u\in GL_n(k)$, which means that all the Hopf algebras associated
to the elements of an orbit of the action \eqref{anactcon}
are isomorphic.
We shall now prove the converse. First we show that two finite
dimensional Hopf algebras $L_1$ and $L_2$ are isomorphic
if and only if $({\cal H}(L_1), {\cal R}_{L_1})$ and
$({\cal H}(L_2), {\cal R}_{L_2})$ are isomorphic as objects
in $\underline{Pent}$. The "only if" part follows from
\thref{te2} and part 5) of \thref{te1}. It remains to prove the
"if" part. Let $f: L_1 \to L_2$ be a Hopf algebra
isomorphism. Then, $f^* :L_2^*\to L_1^*$, $f^*(l^*)=l^*\circ f$ is an
isomorphism of Hopf algebras and
$$
\tilde{f}: {\cal H}(L_1) \to {\cal H}(L_2), \quad
\tilde{f} (h\#h^*):=f(h)\#(f^*)^{-1}(h^*)=f(h)\# h^* \circ f^{-1}
$$
for all $h\in L_1$, $h^*\in L_1^*$ is an algebra isomorphism
between the two Heisenbergs.
On the other hand, if $\{e_i, e_i^*\}$ is a dual basis of
$L_1$, then $\{f(e_i), e_i^* \circ f^{-1}\}$ is a dual basis of
$L_2$ and hence
$(\tilde{f}\ot \tilde{f})({\cal R}_{L_1})={\cal R}_{L_2}$, and this
proves that $\tilde{f}$ is an isomorphism in $\underline{Pent}$.
Let $n_i={\rm dim}(L_i)$, $i=1, 2$. Using \prref{anmorita} we
obtain that
$({\cal H}(L_1), {\cal R}_{L_1}) \cong ({\cal H}(L_2), {\cal R}_{L_2})$
if and only if $({\cal M}_{n_1}(k), R_1)\cong ({\cal M}_{n_2}(k), R_2)$
in $\underline{Pent}$, where $R_i$ is the image of ${\cal R}_{L_i}$
under the algebra isomorphism ${\cal H}(L_i)\cong {\cal M}_{n_i}(k)$.
Now, the two matrix algebras ${\cal M}_{n_1}(k)$ and ${\cal M}_{n_2}(k)$
are isomorphic if and only if $n_1=n_2$ and the Skolem-Noether
theorem tells us that any automorphism $g$ of the matrix algebra
${\cal M}_{n_1}(k)$ is an inner one: there exists $u\in GL_{n_1}(k)$ such
that $g(x)=g_u (x)=uxu^{-1}$. Hence we obtain
that $({\cal M}_{n_1}(k), R_1)\cong ({\cal M}_{n_2}(k), R_2)$ in
$\underline{Pent}$ if and only if $n_1=n_2$ and there exists
$u\in GL_{n_1}(k)$ such that
$R_2=(g_u\ot g_u)(R_1)=(u\ot u)R_1(u\ot u)^{-1}$,
i.e. $R_2$ is equivalent to $R_1$, as needed.
\end{proof}

We shall conclude with a few examples, evidencing first of all
the general method of determining invertible solutions of
the pentagon equation $R\in {\cal M}_n(k)\ot {\cal M}_n(k)$.
Let $(e_{ij})_{i,j=1,n}$ be the canonical basis of ${\cal M}_n(k)$.
An element $R\in {\cal M}_n(k)\ot {\cal M}_n(k)$
can be written as follows
\begin{equation}\eqlabel{mamici}
R=\sum_{i,j=1}^n e_{ij}\ot A_{ij}
\end{equation}
for some matrices $A_{ij}\in {\cal M}_n(k)$. Using the formula
\eqref{AoriB}, $R$ viewed as a matrix in ${\cal M}_{n^2}(k)$,
is given by $R=(A_{ij})_{i,j=1,\cdots, n}$, and we can
quichly check if $R$ is invertible (${\rm det}(R)\neq 0$).
Let us clarify the condition for $R$ to be a solution
of the pentagon equation. Taking into account the
multiplication rule among the elements $(e_{ij})$ we have:
$$
R^{12}R^{13}R^{23}=\sum_{i,j,p,r,s=1}^n
e_{ip}\ot A_{ij} e_{rs}\ot A_{jp} A_{rs}\quad
{\rm and} \quad
R^{23}R^{12}= \sum_{a,b,i,p=1}^n
e_{ip}\ot e_{ab}A_{ip} \ot A_{ab}.
$$
Hence, $R$ is a solution of the pentagon equation if and only if
\begin{equation}\eqlabel{matriciba}
\sum_{a,b=1}^n e_{ab}A_{ip} \ot A_{ab}=
\sum_{i,r,s=1}^n
A_{ij} e_{rs}\ot A_{jp} A_{rs}
\end{equation}
or equivalently,
\begin{equation}\eqlabel{maiscur}
R(A_{ip}\ot I_n)=(\sum_{j=1}^n A_{ij}\ot A_{jp})R
\end{equation}
for all $i$, $p=1, \cdots n$, which can be viewed
as an equation in ${\cal M}_{n^2}(k)$.
We record this observation in the following:

\begin{proposition}\prlabel{recpen}
Let $n$ be a positive integer and
$R=(A_{ij})_{i,j=1,n} \in {\cal M}_{n^2}(k)\cong
{\cal M}_n (k)\ot {\cal M}_n (k)$, an invertible matrix.
Then $R$ is a solution of the pentagon equation if and only if
\begin{equation}\eqlabel{recpenec}
\sum_{j=1}^n A_{ij}\ot A_{jp}=R(A_{ip}\ot I_n)R^{-1}
\end{equation}
for all $i$, $p=1, \cdots, n$.
\end{proposition}

\begin{examples}
\rm
1. Let $n$ be a positive integer,
$A= e_{21} + e_{32} + \cdots + e_{n,n-1} + e_{1n}
\in {\cal M}_n(k)$
and $A_{ij}=\delta_{ij}A^{i-1}$ for all $i$, $j=1,\cdots n$.
Then
$$R=e_{11}\ot I_n + e_{22}\ot A +\cdots + e_{nn} \ot A^{n-1}
\in {\cal M}_n(k)\ot {\cal M}_n(k)$$
is an invertible solution of the pentagon equation
and $P(n, R)\cong (kG)^*$, the Hopf algebra of functions
on a cyclic group with $n$ elements $G$.
Indeed, the pentagon equation \eqref{maiscur} becomes
$$
R(A_{ii}\ot I_n)=(A_{ii}\ot A_{ii})R
$$
for all $i=1,\cdots, n$. As $A_{ii}=A^{i-1}$, it is enough if we
prove that $R(A\ot I_n)=(A\ot A)R$. Using the expression of $A$
and the fact that $A^n=I_n$, we have
$$
R(A\ot I_n)=(A\ot A)R=
e_{1n}\ot I_n + e_{21}\ot A +\cdots + e_{n, n-1}\ot A^{n-1}
$$
i.e. \eqref{maiscur} holds and $R$ is an invertible solution
of the pentagon equation. It remains to prove that
$P(n, R)\cong (kG)^*$. We shall prove that
$H=H(n, R) \cong kG$, the group algebra of $G$ and
then use the duality between $P(n, R)$ and
$H(n, R)$ given by \thref{te1}.
We remark that $R$ is already written in the form
$R=\sum_{i=1}^n A_i \ot B_i$, where $(A_i)$ and $(B_i)$ are
linearly independent. Then $H(n, R)$ is the
commutative subalgebra of ${\cal M}_n(k)$ with
$\{I_n, A, A^2, \cdots, A^{n-1}\}$ as a basis.
Using \eqref{ec919}, \eqref{ec91}
and the fact that
$R^{-1}=e_{11}\ot I_n + e_{22}\ot A^{n-1} +\cdots +
e_{nn} \ot A$,
we obtain that the comultiplication, the counit and the
antipode of $H$ are given by
$$
\Delta_H (A)=A\ot A, \quad \varepsilon_H(A)=I_n, \quad
S_H(A)=A^{n-1}=A^{-1}
$$
i.e. $H\cong kG$.

2. Let $R\in {\cal M}_4(k)\ot {\cal M}_4(k)$ given by
$$
R=(e_{11} + e_{44})\ot I_4 + (e_{22} + e_{33})\ot
(e_{12} + e_{21} + e_{34} + e_{43}) +
e_{13}\ot (e_{31} - e_{42}) + e_{24}\ot (e_{41} - e_{32})
$$
Then $R$ is an invertible solution of the pentagon equation
$H(4, R)\cong H_4$, and hence
$P(4, R)\cong H_4^*\cong H_4$, where $H_4$ is
the Sweedler four dimensional noncommutative noncocommutative
Hopf algebra.

The proof is similar to the previous example
(we left to the reader to check the pentagon
equation \eqref{maiscur} for $R$). The inverse of $R$ is
$$
R^{-1}=(e_{11} + e_{44})\ot I_4 + (e_{22} + e_{33})\ot
(e_{12} + e_{21} + e_{34} + e_{43}) +
e_{13}\ot (e_{41} - e_{32}) + e_{24}\ot (e_{24} - e_{31}).
$$
and therefore the Hopf algebra $H(4, R)$ is the four
dimensional subalgebra of ${\cal M}_4(k)$ having
$$
\{I_4, e_{12} + e_{21} + e_{34} + e_{43},
e_{41} - e_{32}, e_{31} - e_{42} \}
$$
as a $k$-basis. Now, writing
$x=e_{31} - e_{42}$ and $g=e_{12} + e_{21} + e_{34} + e_{43}$
we find that
$$
x^2=0, \quad g^2=I_4, \quad
gx=-xg=e_{41} - e_{32}
$$
On the other hand, the formula of the comultiplication of $H(4,R)$
given by \eqref{ec919}, namely $\Delta (A)= R(A\ot I_4) R^{-1}$ for all
$A\in H(4, R)$, gives, using the expression of $R^{-1}$,
$$
\Delta (g)=g\ot g, \quad \Delta (x)= x\ot g + I_4\ot x
$$
i.e. $H(4, R)\cong H_4$, the Sweedler four dimensional Hopf algebra.
\end{examples}

\section{Conclusions and outlooks}
The study of the pentagon equation, together with a few
classic results in algebra (the Fundamental Theorem for
Hopf modules, Morita theory and the Skolem-Noether theorem), have led
us to the Structure and Classification theorem for finite
dimensional Hopf algebras. \thref{thclas} opens a new road for
describing the types of isomorphisms for Hopf algebras of a certain
dimension. The first step, and the most important, is however the
development of a new Jordan type theory
(we called it {\sl restricted Jordan theory}).
From the point of view of actions,
the classical Jordan theory gives the most elementary description
of the representatives of the orbits of the action
$$
GL_n (k)\times {\cal M}_n(k)\to {\cal M}_n(k),\quad
(U, A)\to UAU^{-1}.
$$
The restricted Jordan theory refers to the following open problem:

{\sl Problem: Describe the orbits of the action}
$$
GL_n(k)\times ({\cal M}_n(k)\ot {\cal M}_n(k))\to
{\cal M}_n(k)\ot {\cal M}_n(k),\quad
(U, R)\to (U\ot U) R (U\ot U)^{-1}.
$$
We recall that the canonical Jordan form $J_A$ of a matrix $A$ is the
matrix equivalent to $A$ which has the greatest number of zeros.
For practical reasons, in the restricted Jordan theory we are in fact
interested in finding the elements of each orbit that have the greatest
number of zeros. Of these, we retain only those which are invertible
solutions of the pentagon equation. The set
of types of $n$-dimensional Hopf algebras shall be those Hopf
algebras associated (using \thref{structurefinite}) to the solutions of
length $n$ (or, equivalently, the space of coinvariant elements is
$n$-dimensional); all other Hopf algebras will have a dimension that is
a divisor of $n^2$. We mention that, as a general rule, the set
of types of $n$-dimensional Hopf algebras is infinite
(this was proved recently in \cite{AndrS}, \cite{BeattieDG},
\cite{Gelaki}). If however we limit ourselves to classifying certain
special types of Hopf algebras, then this set can be finite.
For instance, the set of types of $n$-dimensional semisimple
and cosemisimple Hopf algebras is finite (\cite{Stefan97}).

\end{document}